\documentclass[12pt]{amsart}
\usepackage{amssymb}
\usepackage[T1]{fontenc}
\usepackage[all]{xy}
\UseRawInputEncoding
\numberwithin{equation}{section}

\newcommand{\Q}{{\mathbb {Q}}}

\newcommand{\R}{{\mathbb{R}}}
\newcommand{\Z}{{\mathbb{Z}}}
\newcommand{\C}{{\mathbb{C}}}

\newcommand{\N}{{\mathbb{N}}}

\theoremstyle{plain}
\newtheorem{thm}{Theorem}[section]

\newtheorem{prop}[thm]{Proposition}
\newtheorem{cor}[thm]{Corollary}

\theoremstyle{definition}
\newtheorem{remark}[thm]{Remark}

\title[Norm forms]{Characterization of norm forms via their values at integer points}

\author{George Tomanov}
\address{Institut Camille Jordan, Universit\'e Claude Bernard - Lyon
I, B\^atiment de Math\'ematiques, 43, Bld.
 du 11 Novembre 1918,
69622 Villeurbanne Cedex, France {\tt tomanov@math.univ-lyon1.fr}}

\begin{document}

\maketitle

\begin{abstract} Using a homogeneous dynamics approach, we obtain a complete description of the forms with discrete set of values at
 integer points and not representing zero over the rational numbers non-trivially. It turns out that these are exactly the norm forms and the quasi-norm forms (introduced in the paper).
 As a by-product, we obtain a general class of non-totaly real forms for which a natural generalization of the conjecture of Cassels and Swinnerton-Dyer fails.

\end{abstract}

\section{Introduction} \label{Introduction}

The main purpose of this paper is to give a characterization of
the norm forms and related quasi-norm forms (introduced in this paper) in terms of their values at
 integer points. Note that the norm forms represent a classical object of study in both
 algebraic number theory \cite[ch.2]{BoSha} and Diophantine approximation \cite[ch.7]{Schmidt}. From a contemporary viewpoint, the investigation of similar forms
is a topic of the dynamics of group actions on homogeneous spaces. (See \cite{Margulis problems survey}, \cite{Oppenheim survey}, \cite{EL}, \cite{Toma3}.)

Let us introduce our basic notations.
We denote by $\R[\vec{x}]$ the ring of real polynomials in $m \geq 2$ variables $\vec{x} = (x_1, \dots, x_m)$.
In what follows, by a \textit{form} we mean a homogeneous polynomial $f(\vec{x}) \in \R[\vec{x}]$ which factors as $f(\vec{x})= l_1(\vec{x}) \dots l_n(\vec{x})$, where $l_1(\vec{x}), \dots,
l_n(\vec{x})$ are linearly independent \textit{linear} forms with complex coefficients.
 It is clear that $n = 2s + t$, where $s$ is the number of pairs of complex-conjugate linear forms, and $t$ is the number of real linear forms. (\textit{Further on, given a form $f$ of degree $n$, the notations $s$ and $t$ have the same meaning}.) The form $f$ is called \textit{totally real} if $s=0$, and \textit{totally imaginary} if $t=0$. Let $G = \mathrm{SL}_m(\R)$.  The group $G$ is considered as a (real) algebraic group defined over $\Q$, shortly, a $\Q$-group. So, $G(\Q) = \mathrm{SL}_m(\Q)$, and $\Gamma = G(\Z)$ is a lattice in $G$.
  The  homogeneous space
$G/\Gamma$ is endowed with the quotient topology. As usual, $\pi: G \rightarrow G/\Gamma, g \mapsto g\Gamma$, is the natural projection. The group $G$ is acting on $G/\Gamma$ by left translations, and $G$ is acting on $\R[\vec{x}]$ by $(g p)(\vec{x}) = p(g^{-1}\vec{x})$, where $g \in G$ and $p(\vec{x}) \in \R[\vec{x}]$. We denote by $H$ the Zariski identity component of the stabilizer of the form $f$ in $G$.
Note that $m \geq n$. If $m = n$, then $H$ is a maximal algebraic torus of $G$ and $H$ is an almost direct product of its maximal compact sub-torus $H_c$ and of its maximal
$\R$-\textit{split} (or \textit{diagonalizable over} $\R$) sub-torus $H_d$. (See subsection \ref{2.3} for details.)
Observe that a form in $n$ variables $f$ is totally real if and only if $H = H_d$, in which case $H$ is conjugate to the diagonal subgroup of $G$.
If $g \in G(\Q)$, then the forms $f$ and $g f$ are called $\Q$-\textit{equivalent}. If, in addition,
$(gf)(x_1, \dots, x_m) = (g f)(x_1, \dots, x_n, 0, \dots, 0)$, we say that $f$ is $\Q$-equivalent to a form in $n$ variables.
Further on, we denote by $\mathcal{Z}$
a $\Gamma$-invariant subset of $\Z^m \setminus \{\vec{0}\}$. It is clear that $\mathcal{Z}$ is a disjoint union of \textit{level sets} $\mathcal{Z}_d, d \in \N$, defined by $\mathcal{Z}_d = d \mathcal{Z}_1$, where $\mathcal{Z}_1 = \{(z_1, \dots, z_m) \in \Z^m: \mathrm{gcd}(z_1, \dots, z_m) = 1\}$ is the set of primitive integral points.
 We are interested in the study of forms with discrete $f(\mathcal{Z})$.
  Our first result reduces the study (up to $\Q$-equivalence) of the forms $f$ to the case when $m = n$.

\begin{thm}
\label{thm0} If $f(\mathcal{Z})$ is discrete in $\R$, then $f$ is $\Q$-equivalent to a form in $n$ variables.
\end{thm}

\textit{Up to the end of the Introduction, we suppose that } $m = n$. 
If $f$ is totally real and $f(\mathcal{Z})$ is discrete, then, in view of \cite[Theorem 1.8]{Toma2}, $f$ is a constant multiple of a form with rational coefficients.
The case when $f$ is not totally real is essentially different, and the  analogous result is not valid.
For example, if $f_0(\vec{x}) = q_1(x_1,x_2) \cdots q_s(x_{2s-1},x_{2s})x_{2s+1}\cdots x_{n}$, where $q_i$ are definite quadratic forms in two variables, $n \geq 3$, $s \geq 1$, and at least one of the forms $q_i$ is not proportional to a form with rational coefficients, then
$f_0$ is also not proportional to a form with rational coefficients despite the discreteness of  $f_0(\mathcal{Z})$.
Observe that $f_0$
represents $0$ over $\Q$ non-trivially, that is, $f(\vec{v}) = 0$ for some $\vec{v} \in \Q^n \setminus \{\vec{0}\}$.
It turns out that the behaviour of $f(\mathcal{Z})$ around $0$ is related to important "rigidity" properties separating the totally real from the non-totally real forms and giving rise to notable conjectures. We say that $f(\mathcal{Z})$ \textit{is discrete at} $0$ if $f(\mathcal{Z}) \cap (-\varepsilon, \varepsilon) = \{0\}$  for some $\varepsilon > 0$.

In the present context, it is worthwhile to mention the 1955 conjecture of Cassels and Swinnerton-Dyer, implying the famous conjecture of Littlewood  stated around 1930.  (See \cite{CS}, \cite[Conjectures 7 and 8]{Margulis problems survey} or \cite[\S2]{Oppenheim survey}.)

\medskip

\textbf{Conjecture 1}. \textit{Let $f$ be a totally real form such that
\begin{enumerate}
  \item[(i)] $n \geq 3$;
\item[(ii)]
$f$ does not represent $0$ over $\Q$ non-trivially;
\item[(iii)] $f(\Z^n)$ is discrete at $0$.
\end{enumerate}
Then $f$ is proportional to a form with rational coefficients.}

\medskip

An example of a form in two variables justifying the assumption $n \geq 3$ in the formulation of the conjecture
is given by $\alpha^2 x_1^2 - x_2^2$ , where $\alpha$ is a badly approximable number such that $\alpha^2 \notin \Q$.
(See \cite[1.2]{Margulis problems survey} and \cite[ch.1, \S5]{Schmidt}.)

Conjecture 1 can be reformulated in terms of homogeneous dynamics as follows \cite[Conjecture 9]{Margulis problems survey}:

\medskip

\textbf{Conjecture 2}.
\textit{If $n \geq 3$ and $H$ is a maximal diagonalizable over $\R$ algebraic torus of $G$, then every relatively compact orbit $H\pi(g)$ is compact.}

\medskip

In contrast to Conjecture 1, \textit{in this paper we treat general (not necessarily totally real) forms} $f$.
We describe the forms $f$  not representing $0$ over $\Q$ non-trivially and with discrete $f(\mathcal{Z})$.
As a by-product, we obtain a general construction of non-totally real forms $f$ in $n = 2s + t \geq 3$ variables, with no additional restrictions on $s$ and $t$, for which a natural generalization of Conjecture 1 does not hold.

In order to formulate our main result, we recall the notion of a norm form \cite[ch.7, \S1]{Schmidt} and introduce the notion of a quasi-norm form.
Let $K = \Q \alpha_1 + \dots + \Q \alpha_n$ be a number field of degree $n$ and $\{\theta_1, \dots, \theta_n\}$
be the set of all embeddings of $K$ into $\C$.
Denote
$$l_i(\vec{x}) = \theta_i(\alpha_1)x_1 + \dots + \theta_i(\alpha_n)x_n.$$ 
Let
$\mathfrak{N}_{K}(\vec{x}) = \overset{n}{\underset{i=1}{\prod}}l_{i}(\vec{x})$. It is clear that $\mathfrak{N}_{K}(\vec{x}) \in \Q[\vec{x}]$.
By a \textit{norm form} \textit{corresponding to} $K$, we mean a form $\Q$-equivalent to the following one:
 \begin{equation}
\label{def1}
f(\vec{x}) = c\cdot \mathfrak{N}_{K}(\vec{x}), \  c \in \R^*.
\end{equation}

The quasi-norm forms are totally imaginary forms, defined as follows.
Let $F$ be a totally real number field of degree $s \geq 1$, $\{\theta_1, \dots, \theta_{s}\}$ be the set of all real embeddings of $F$ and $F = \Q\alpha_1 + \dots + \Q\alpha_{s}$. Denote $\lambda_i(\vec{y}) = \theta_i(\alpha_1)y_1 + \dots + \theta_i(\alpha_{s})y_{s}$, $1 \leq i \leq s$, where $\vec{y} = (y_1, \dots, y_s)$.
Let $q_1, \dots, q_{s}$ be definite quadratic forms on $\R^2$. Put $n = 2s$ and
write $\vec{x} = (\vec{x}_1, \vec{x}_2)$, where $\vec{x}_1 = (x_1, \dots, x_{s})$ and $\vec{x}_2 = (x_{s+1}, \dots, x_{n})$.
By a $\textit{quasi-norm form corresponding to}$ $F$ we mean a form $\Q$-equivalent to
\begin{equation}\label{def2}
  f(\vec{x}) = q_1(\lambda_1(\vec{x}_1), \lambda_1(\vec{x}_2)) \cdots q_s(\lambda_{s}(\vec{x}_1), \lambda_{s}(\vec{x}_2)).
\end{equation}

The classes of norm forms and quasi-norm forms have a tiny intersection.
As elucidated by Remark \ref{form=quasi/form}, $f$ is both a norm and a quasi-norm form if and only if $f$ is given (up to $\Q$-equivalence) by
(\ref{def1}), where $K$ is a \textit{complex multiplication number field} (shortly, a $\mathrm{\textit{CM}}$ \textit{field}), that is, $K = F(\sqrt{-a})$, where $F$ is  a totally real number field and $a \in F$ is such that $\theta_i(a) > 0$ for all $\theta_i$. In this case,
$f$ is also
given (up to $\Q$-equivalence) by (\ref{def2}) with $q_i(y_1,y_2)$ proportional to $y_1^2 + \theta_i(a)y_2^2$.

We now state the main result of this paper.

\begin{thm}
\label{thm++} The following conditions are equivalent:
\begin{enumerate}
  \item $f$ is either a norm form or a quasi-norm form;
  \item $f(\mathcal{Z})$ is discrete and $f$ does not represent $0$ over $\Q$ non-trivially;
 \item $H\pi(e)$ is compact. In this case, we have in terms of $H$:
 \begin{enumerate}
  \item $f$ is a norm form  if and only if 
  $H \cap G(\Q)$ is dense is $H$;
  \item $f$ is a quasi-norm form if and only if $n = 2(\dim H_d + 1)$ and $H_d\pi(e)$ is compact.
\end{enumerate}
\end{enumerate}
\end{thm}

It is worth mentioning that, besides some similarities,  the proofs of Theorems \ref{thm0} and \ref{thm++} are basically different from the aforementioned \cite[Theorem 1.8]{Toma2}. For instance, in the present paper, we do not need a description of divergent trajectories of diagonal subgroups in $G/\Gamma$ (see \cite{TW}). Instead, we essentially use
Mahler's compactness criterion, theorems for orbit closures, and results and techniques from the theory of artinian rings.

Since every quasi-norm form is totally imaginary, Theorem \ref{thm++} immediately implies the following.
\begin{cor}
\label{cor1} Let $f$ be a form such that
\begin{enumerate}
  \item[(i)]
$f(\mathcal{Z})$ is discrete;
\item[(ii)]
$f$ does not represent $0$ over $\Q$ non-trivially;
\item[(iii)] $f$ is not totally imaginary, which, in particular, is satisfied if $n$ is odd.
\end{enumerate}
Then $f$ is a norm form.
\end{cor}

Corollary \ref{cor1} shows that the totally real form $f$ in the formulation of Conjecture 1 is a norm form. Since we are interested in arbitrary forms, \textit{assume that $f$ is not totally real (that is, $s > 0$) and satisfies conditions $(i)-(iii)$, as in the formulation of Conjecture 1}. In view of Theorem \ref{thm++}, the natural question arises as to whether $f$ is either a norm or a quasi-norm form provided some additional hypothetical restrictions on $s$ and $t$ are also satisfied. In subsection $5.1$ we introduce a general class of non-totally real forms giving a negative answer. In more concrete terms, we get:

\begin{thm} With $n \geq 3$, the following holds.
\label{cor2}
\begin{enumerate}
  \item
 Let $n = 2s + t$ with $s > 0$. There exists $f$ of degree $n$, with $s$ the number of pairs of complex-conjugate linear forms dividing $f$, which is neither a norm nor a quasi-norm form, although $f$ does not represent $0$ over $\Q$ non-trivially and $f(\mathcal{Z})$ is discrete at $0$;
\item If $H$ is a maximal non-diagonalizable over $\R$ algebraic torus of $G$, there exists an orbit  ${H\pi(g)}$, which is relatively compact but not compact.
\end{enumerate}
\end{thm}

Note that the analogous statements are not valid for $n = 2$ because, in this case, every non-totally real form $f$ is a quasi-norm form corresponding to $F = \Q$, and every maximal non-diagonalizable over $\R$ algebraic torus $H$ is compact.

The structure of the paper is the following. Section 2 assembles auxiliary results from the theory of algebraic $\Q$-groups  and their arithmetic subgroups and from the algebraic number theory. 
Theorem \ref{thm0} is proved in Section 3 using  Ratner's result \cite{Ratner} and its specification \cite[Theorem 1]{Toma1}. Theorem \ref{thm++} is proved in Section 4, where we establish compactness criteria for orbits of tori on $G/\Gamma$, reflecting the structure of $f$ and of its stabilizer in $G$ (Propositions \ref{auxiliary2} and \ref{anisotropic1}).
 Using Theorem \ref{thm++}, Theorem \ref{cor2} is proved in Section 5.

\section{Preliminaries}

\subsection{Forms and their stabilizers.}\label{2.1}
Any vector of $\R^m$ is identified with its coordinate vector with respect to the canonical basis $\vec{e}_1, \dots, \vec{e}_m$ of the $\Q$-vector space $\R^m$. If $g \in G$ and $\vec{x} = (x_1, \dots, x_m) \in \R^m$, then $g \vec{x}$ is given by the matrix product $\vec{x}\cdot g^\tau$, where $g^\tau$ is the transpose of the matrix $g$.
The group $G = \mathrm{SL}_m(\R)$ is considered as a real $\Q$-algebraic group, that is, as the group of $\R$-points of the algebraic group defined by the polynomial $\det (x_{ij})_{1 \leq i,j \leq m} - 1$. We refer to \cite{Borel} for the basic notions and facts about algebraic groups. All topological notions connected with the Zariski topology will be preceded by the prefix "Zariski".

As in the introduction, $f(\vec{x}) \in \R[\vec{x}]$ is a decomposable  over $\C$ homogeneous form in $m \geq 2$ variables $\vec{x} = (x_1, \dots, x_m)$, namely,
$f(\vec{x})= l_1(\vec{x}) \cdots l_n(\vec{x})$, where $l_1(\vec{x}), \dots,
l_n(\vec{x})$ are linearly independent over $\C$ linear forms with complex coefficients. 
We assume that $l_{2i}(\vec{x}) = \overline{l}_{2i-1}(\vec{x})$, $1 \leq i \leq s$, where $\overline{l}_{2i-1}$ is the linear form conjugate to ${l}_{2i-1}$, and that $l_{2s+j} \in \R[\vec{x}], 1 \leq j \leq t$. (Recall that $n = 2s + t$.) It is easy to see that there exist
linearly independent real linear forms $\{\lambda_1, \dots, \lambda_n \}$ such that ${l}_{2i-1}\cdot \overline{l}_{2i} = \lambda_{2i-1}^2 + \lambda_{2i}^2, 1 \leq i \leq s$, and $\lambda_{2s+j} = l_{2s+j}, 1 \leq j \leq t$.
Replacing, if necessary, $f$ by $c\cdot f, c \in \R^*$, we can fix
a $\sigma \in G$ such that $\sigma^{-1}\vec{x} = (\lambda_1, \dots, \lambda_n, x_{n+1}, \dots, x_m)$.
Then
$f(\vec{x}) =  (\sigma f_0)(\vec{x})$, where
\begin{equation}
\label{form}
f_0(\vec{x})= (x_{1}^2 + x_{2}^2) \cdots (x_{2s-1}^2 + x_{2s}^2)  x_{2s + 1} \cdots x_n.
\end{equation}

Let $H_0$ be the Zariski identity component of the stabilizer of $f_0$ in $G$.
Then
\begin{equation}
\label{stab}
H_0 = \left\{\left(\begin{array}{cc} a&0\\ b&c\\ \end{array} \right):a \in S, \ b \in \mathrm{M}_{m-n,n}(\R), \ c \in \textrm{SL}_{m-n}(\R) \right\},
\end{equation}
where $S$ is the Zariski identity component of the subgroup of $\mathrm{SL}_n(\R)$ fixing the restriction of $f_0$ on $\R^n$.
Note that $H = \sigma H_0 \sigma^{-1}$ is the Zariski identity component of the stabiliser of $f$.
\subsection{Maximal $\R$-tori of $G$.}\label{2.3} 
Up to the end of this subsection, we suppose that $m = n \geq 3$ and use the notation of subsection 2.1.
In this paper, by a \textit{torus} we mean \textit{a torus of the real algebraic group} $G$, that is, a Zariski closed and Zariski connected subgroup of $G$, which is diagonalizable over $\C$. If $L$ is a subgroup of $G$ closed in the Hausdorff topology, we denote by $L^\circ$ its connected component of the identity.

Let
$$\R^n = V_{0,1} \oplus \dots \oplus V_{0,s+t}$$
be the direct sum of irreducible $H_0$-submodules of $\R^n$. 
In view of (\ref{form}), $V_{0,i} = \R e_{2i-1} + \R e_{2i}$ if $1 \leq i \leq s$, and
$V_{0,i} = \R e_{s + i}$ if $s <  i \leq s+t$. The groups $\mathrm{GL}(V_{0,i})$ are considered as subgroups of $\mathrm{GL}(\R^n)$, and the rings $\mathrm{End}(V_{0,i})$ are considered as subrings of the ring of all $n \times n$ real matrices $\mathrm{M}_n(\R)$ identified with $\mathrm{End}(\R^n)$. (Below we implicitly use the inclusion $G \subset \mathrm{GL}(\R^n)$.)
For every $1 \leq i \leq s + t$ denote by  $Z_{0,i}$ the center of $\mathrm{GL}(V_{0,i})$ (isomorphic to $\R^*$), and
for every $1 \leq i \leq s$ denote by $q_{0,i}$ the restriction of the quadratic form $x_{2i - 1}^2 + x_{2i}^2$ to $V_{0,i}$.
Let $H_{0,i} = \mathrm{SO}(q_{0,i}) \cdot Z_{0,i}$ if $1 \leq i \leq s$, and $H_{0,i} = Z_{0,i}$ if $s < i \leq s + t$.
Obviously, $H_{0,i} \cong \C^*$ if $1 \leq i \leq s$, and $H_{0,i} \cong \R^*$ if $s < i \leq s+t$.
Since $(H_{0,1} \times \dots \times H_{0,s+t}) \cap G$ fixes $f_0$ and contains $H_0$, we have
$$
H_0 = (H_{0,1} \times \dots \times H_{0,s+t}) \cap G.
$$
Clearly $H_0$ (as well as $H = \sigma H_0 \sigma^{-1}$) is a maximal torus of $G$.

Let $\mathcal{M}_0$ be the subspace of the vector space $\mathrm{M}_n(\R)$ generated by $H_0$. Then
$$
\mathcal{M}_0 = \mathcal{M}_{0,1} \oplus \dots \oplus \mathcal{M}_{0,s+t},
$$
where $\mathcal{M}_{0,i} \subset \mathrm{End}(V_{0,i})$, $\mathcal{M}_{0,i}\cong \C$ for $1 \leq i \leq s$, and $\mathcal{M}_{0,i} \cong \R$ for $s < i \leq s+t$.

 Denote $H_{i} = \sigma H_{0,i}\sigma^{-1}$, $Z_{i} = \sigma Z_{0,i}\sigma^{-1} (\cong \R^*)$, $V_{i} = \sigma V_{0,i}$, $\mathcal{M}_{i} = \sigma \mathcal{M}_{0,i}\sigma^{-1}$, $i = 1, \dots, s+t$, and $\mathcal{M} = \sigma \mathcal{M}_{0}\sigma^{-1}$. 
 Let $q_i$ be the restriction of $\lambda_{2i-1}^2 + \lambda_{2i}^2$ to $V_i$, $1 \leq i \leq s$. Then $q_i$ is a positive definite quadratic form on $V_i$, and $H_i = \mathrm{SO}(q_{i}) \cdot Z_{i} \cong \C^*$, $1 \leq i \leq s$. Also, $\mathcal{M}_{i} \cong \C$ when $1 \leq i \leq s$, and $\mathcal{M}_{i} \cong \R$ when $i > s$. Therefore
$$
H = (H_{1} \times \dots \times H_{s+t}) \cap G \subset \mathcal{M} = \mathcal{M}_{1} \oplus \dots \oplus \mathcal{M}_{s+t} \cong \C^s \oplus \R^{t}.
$$
\textit{We will denote by} $p_i$ \textit{the projection} $\mathcal{M} \rightarrow \mathcal{M}_{i}$ \textit{parallel to} $\underset{j \neq i}{\oplus}\mathcal{M}_{j}$.

Let $H_c = \mathrm{SO}(q_{1}) \times \dots \times \mathrm{SO}(q_{s})$ and $H_d = (Z_1 \times \dots \times Z_{s + t}) \cap G$. Then
$H_c$ is the maximal compact sub-torus of $H$, $H_d$ is the maximal $\R$-split sub-torus of $H$, and $H$ is an almost direct product of $H_c$ and $H_d$, that is, $H = H_c \cdot H_d$ and $H_c \cap H_d$ is finite.
It is easy to see that the tori $H$ (defined in this subsection) exhaust all maximal tori of the real algebraic group $G$.
The $\R$-rank of $H$ is equal to $s+t-1$ and $H$ is determined by its $\R$-rank up to conjugation in $G$.

Let $\R[\vec{x}]_n$ be the subspace of $\R[\vec{x}]$ generated by all homogeneous polynomials of degree $n$, $H$ be a maximal torus of $G$, and $\R[\vec{x}]_n^{H}$ be the subspace of $\R[\vec{x}]_n$ consisting of all $H$-invariant polynomials.
There exists a form $f$ such that $\R[\vec{x}]_n^{H} \supset \R f$. The following proposition establishes a bijective correspondence between the maximal tori $H$ of $G$ and the classes $\R f$ of proportional forms.
\begin{prop}
\label{prop.inv}
With the above notation,
\begin{equation}
\label{i-form}
\R[\vec{x}]_n^{H} = \R f.
\end{equation}
\end{prop}
{\bf Proof.} Let $D$ be the diagonal subgroup of $\mathrm{SL}_n(\C)$, $\C[\vec{x}]_n = \R[\vec{x}]_n \otimes_{\R} \C$ and $h_0 = x_1 \cdots x_n$. Let $\C[\vec{x}]_n^D$ be the subspace of all $D$-invariant polynomials of $\C[\vec{x}]_n$. It is easy to see that
$\C[\vec{x}]_n^D = \C h_0$. There exists $g \in \mathrm{SL}_n(\C)$ such that $gh_0 = a\cdot f$, $a \in \C^*$. Hence $gDg^{-1}$ is the Zariski closure of $H$ in $\mathrm{SL}_n(\C)$ which implies that
$$
\C f = \C[\vec{x}]_n^{gDg^{-1}} = \C[\vec{x}]_n^{H}.
$$
Since $\C[\vec{x}]_n^{H} = \C \otimes_{\R} \R[\vec{x}]_n^{H}$, we get (\ref{i-form}). 
\qed

\subsection{$\Q$-tori of $G$.}\label{2.4}
A torus $T$ of $G$ (not necessary maximal) is \textit{defined over} $\Q$ (or is a $\Q$-\textit{torus}) if its subgroup of $\Q$-points $T(\Q) = T \cap G(\Q)$ is dense in $T$ \cite[18.3]{Borel}. Also, a $\Q$-torus $T$ is $\Q$-\textit{anisotropic} if it does not admit non-trivial $\Q$-rational characters; equivalently, if $T(\Z) = T \cap \Gamma$ is co-compact in $T$.
A $\Q$-torus $T$ is $\Q$-\textit{split} if it is diagonalizable over $\Q$; equivalently, if every character of $T$ is defined over $\Q$ \cite[8.2]{Borel}.
Every $\Q$-torus is an almost direct product of its maximal $\Q$-anisotropic sub-torus and of its maximal $\Q$-split sub-torus \cite[8.15]{Borel}.

\subsection{Fields with "units defect".}\label{2.4+}
The following proposition is well-known. For reader's convenience, we provide a short proof of it.

\begin{prop}
\label{fields}
Every maximal subfield $K$ of ${M}_n(\Q)$ is a number field of degree $n$.
\end{prop}
{\bf Proof.}
By the theorem of double centralizer \cite[Theorem 4.3.2]{Herstein}, the centralizer  $C_{\mathrm{M}_n(\Q)}(K)$ of $K$ in $\mathrm{M}_n(\Q)$ is isomorphic to $\mathrm{M}_r(\mathcal{D})$, where $\mathcal{D}$ is a division algebra with center $K$. Since $K$ is maximal, $\mathcal{D} = K$ and $r = 1$, that is, $C_{\mathrm{M}_n(\Q)}(K) = K$. Hence, $C_{\mathrm{M}_n(\Q)\otimes_\Q \C}(K \otimes_\Q \C) = K \otimes_\Q \C \cong \C^n$, proving that $\dim_{\Q} K = n$. \qed

\medskip

\textit{Further on, we will use the notation}
\begin{equation}
\label{Delta}
\Delta = H \cap \Gamma.
\end{equation}

Let $K$ be a maximal subfield of ${M}_n(\Q)$ and $f = \mathfrak{N}_K$. Then $f$ is the restriction of $\det$ to
$K$. Let $H$ be the maximal torus of $G$ corresponding to $f$ (Proposition \ref{prop.inv}). Then $H$ is a  maximal $\Q$-torus of $G$. With the notation from subsection $2.2$, we have $K \subset \mathcal{M}$ and $H(\Q) = K \cap G$. The restrictions of $p_i$, $1 \leq i \leq s + t$, to $K$ are the complex embeddings and the restrictions of $p_i$, $s < i \leq s + t$, to $K$ are the real embeddings of $K$. Recall that $K$ is called \textit{totally imaginary} (respectively, \textit{totally real}) if $t = 0$ (respectively, $s = 0$). If
$\mathcal{O}_K$ is the ring of integers of $K$ and $\mathcal{O}_K^*$ is the unit group of $\mathcal{O}_K$, then $K \cap \mathrm{M}_n(\Z)$ is a subgroup of finite index in $\mathcal{O}_K$ and $\Delta$ is a subgroup of finite index in $\mathcal{O}_K^*$.
By Dirichlet's unit theorem \cite[Ch.2, Theorem 5]{BoSha}, $\Delta$ is of rank $s+t-1$, that is, $\Delta$ contains a  subgroup of finite index isomorphic to $\Z^{s+t-1}$. We say that the field $K$ has "units defect" if it contains a \textit{proper} subfield such that its unit group has rank $s+t-1$.

Recall the following result from \cite{Remak}.

\begin{prop}
\label{CM-field}
$K$ has "units defect" if and only if $K$ is a $CM$ field.
\end{prop}

Note that $H_i = \mathrm{SO}(q_{i}) \times {Z}_{i}^\circ, 1 \leq i \leq s,$ where $\mathrm{SO}(q_{i}) \cong \mathrm{SO}_2(\R)$
and ${Z}_{i}^\circ \cong \R_{>0}$. Let $r_i: H_i \rightarrow \mathrm{SO}(q_{i})$ be the natural projection.

\begin{prop}
\label{ort.proj.} Assume that
$K$ is not a $CM$ field. Then $r_i \circ p_i (\Delta)$ is dense in $\mathrm{SO}(q_{i})$ for every $1 \leq i \leq s$.
\end{prop}
{\bf Proof.} Suppose to the contrary that there exists $1 \leq i \leq s$ such that $r_i \circ p_i (\Delta)$ is not dense in $\mathrm{SO}(q_{i})$. Then $r_i \circ p_i (\Delta)$ is finite, and, therefore, there exists an integer $a > 0$ such that $p_i(\Delta^{a}) \subset \{e\} \times {Z}_{i}^\circ$, where $\Delta^{a} = \{\xi^a : \xi \in \Delta\}$.
Let $K'$ be the subfield of $K$ generated by $\Delta^{a}$. Then $p_i(K')$ is contained in the subfield of real numbers of $\mathcal{M}_i$ $(\cong \C)$. By the weak approximation theorem \cite[Lemma 6.1]{Cassels}, $p_i(K)$ is dense in $\mathcal{M}_i$, which implies that $K' \subsetneqq K$. It follows from the definition of $K'$ that the unit groups of $K'$ and $K$ are of the same rank. By Proposition \ref{CM-field} $K$ is a $CM$ field, giving a contradiction.\qed

The following general proposition is relevant to Theorem \ref{cor2}.
\begin{prop}
\label{number fields}
Let $n = 2s+t \geq 2$, where $n, s$, and $t$ are non-negative integers. There exists a number field $K$ of degree $n$ with $2s$ complex (non-real) embeddings and $t$ real embeddings.
Moreover, if $t = 0$, then $K$ can be chosen to be not a $CM$ field.
\end{prop}

{\bf Proof.}  The case $n = 2$ being trivial, further on we suppose that $n \geq 3$.
Let $p_{\infty}(x) = x^n + a_1 x^{n-1} + \dots +  a_n \in \R [x]$
be a separable polynomial with $2s$ complex and $t$ real roots. Let $l$ be a prime number and $p_l (x) = x^n + b_1x^{n-1} + \dots + b_n \in \Z_l[x]$, where $\Z_l$ is the ring of integers of the field of $l$-adic numbers $\Q_l$. We suppose that $l \mid b_i$ for all $i$ and $l^2 \nmid  b_n$. So, $p_l (x)$ is irreducible in $\Q_l[x]$ by Eisenstein's criterion. Given $0 < \epsilon \leq \frac{1}{l^2}$, it follows from the weak approximation \cite[Lemma 6.1]{Cassels} that there exists a
$$p(x) = x^n + c_1 x^{n-1} + \dots +  c_n \in \Q [x]$$
such that $\mid c_i - a_i \mid_{\infty} < \epsilon$ and $\mid c_i - b_i \mid_{l} < \epsilon$. (As usual, $\mid \cdot \mid_{\infty}$ is the  euclidean norm on $\R$, and $\mid \cdot \mid_l$ is the $l$-adic norm on $\Q_l$.) By Eisenstein's criterion again, $p(x)$ is irreducible in $\Q [x]$. Let $n = 2s' + t'$, where $2s'$ is the number of complex roots, and $t'$ is the number of real roots of $p(x)$. Fix reals $\alpha_0 < \alpha_1 < \dots < \alpha_t$ such that every interval $(\alpha_{i-1}, \alpha_i)$, $0 < i \leq t$, contains a real root of $p_{\infty}$. So, $p_{\infty}(\alpha_{i-1}) p_{\infty}(\alpha_{i}) < 0$, and, choosing $\epsilon$ small enough, we get that $p(\alpha_{i-1}) p(\alpha_{i}) < 0$. Therefore  $t' \geq t$, proving the proposition for $s = 0$. Further on, assume $s > 0$.
Let $\Phi: \C^n \rightarrow \C^n, \vec{z} \mapsto (\phi_1(\vec{z}), \dots, \phi_n(\vec{z}))$, where $\phi_i$ is the elementary symmetric polynomial in $n$ variables of degree $i$. Let $\vec{\theta} = (\theta_1, \dots, \theta_n)$, where $\theta_i$
are the roots of $p_{\infty}$. Then $\Phi(\vec{\theta}) = \vec{a}$, where $\vec{a} = (a_1, \dots, a_n)$, and $\Phi^{-1}(\vec{a})$
is a finite set of $n!$ vectors, each of them with coordinates obtained by a permutation of the coordinates of $\vec{\theta}$. It follows easily from the definition of $\Phi$ that there exists a neighbourhood $W$ of $\vec{\theta}$ such that $\Phi(W)$ is a neighbourhood of $\vec{a}$ and $\Phi|_{W}: W \rightarrow \Phi(W)$ is a homeomorphism. Choosing $W$ small enough, we get that if $(z_1, \dots, z_n) \in W$ and $\theta_i \notin \R$, then $z_i \notin \R$. Let $\epsilon$ be so close to $0$ that $\vec{c} = (c_1, \dots, c_n) \in \Phi(W)$.
Then if $\vec{\xi} = (\xi_1, \dots, \xi_n) \in W$ is such that $\Phi(\vec{\xi}) = \vec{c}$, the numbers $\xi_i$ are the roots of $p$, and, in view of the choice of $W$,  $K = \Q(\xi_1)$ is a field with exactly $2s$ complex (non-real) embeddings.

It remains to prove that if $n = 2s \geq 4$, there exists a totally imaginary non-$CM$ field of degree $n$.
Let $K = \Q(\theta)$, where $\theta$ is a root of $x^4 + x^3 - x^2 - x + 1 \in \Q[x]$. It is known (and can be proved as an exercise) that $K$ is a totally imaginary non-$CM$ field of degree $4$.
Let
$s \geq 3$ and $s = 2s_1 + t_1$ with $s_1 > 0$ and $t_1 > 0$.
By the first part of the proposition, there exists a number field $K_1$ of degree $s$ with $2s_1$ complex embeddings and $t_1$ real embeddings. Let $\sigma_i: K_1 \rightarrow \R$, $1 \leq i \leq t_1$, be the real embeddings. By the weak approximation, there exists $\beta \in K_1$ such that $\sigma_i(\beta) < 0$ for all $i$. Clearly, $K = K_1(\sqrt{\beta})$ is as needed. \qed

\section{Proof of Theorem \ref{thm0}}

\subsection{On the stabilizer of a polynomial} \label{section:
general proposition}

The following proposition, representing some independent interest, is proved in the particular case of totally real forms $f$ in
 \cite[Proposition 7.3]{Toma2}. The present version is valid for \textit{any} polynomial, its proof is simpler than \cite[Proposition 7.3]{Toma2}, and it remains applicable in the general $S$-adic setting.

\begin{prop}
\label{connection1} Let $\mathrm{p} \in \R[\vec{x}]$, and let $L$ be the stabilizer of $\mathrm{p}$ in $G$. Suppose that $\mathrm{p}(\mathcal{Z})$ is
discrete in $\R$. Then $L\pi(e)$ is closed in
$G/\Gamma$.
\end{prop}

{\bf Proof.} Let $h_i \in L$ be a sequence such that $\lim_{i\rightarrow \infty} h_i\pi(e) = \pi(g)$. We need to prove that $\pi(g) \in L\pi(e)$. There exists a sequence $\tau_i \in G$ converging to $e$ such that
\begin{equation}\label{1}
 h_i\pi(e) = \tau_i \pi(g);
\end{equation}
equivalently,
$$
h_i \gamma_i = \tau_i g,
$$
where $\gamma_i \in \Gamma$. Let $\vec{z} \in \mathcal{Z}$. Then
$$
\mathrm{p}(\gamma_i \vec{z}) = \mathrm{p}(\tau_i g\vec{z})
$$
and, therefore, $\mathrm{p}(\gamma_i\vec{z})$ converges to $\mathrm{p}(g\vec{z})$. Since $\mathcal{Z}$ is $\Gamma$-invariant and $\mathrm{p}(\mathcal{Z})$ is discrete in $\R$, there is a $c(\vec{z})
> 0$ such that
\begin{equation}
\label{2} \mathrm{p}(\tau_i g\vec{z}) = \mathrm{p}(g\vec{z})
\end{equation}
for all $i > c(\vec{z})$.

Let $r = \deg \mathrm{p}$, and let $V$ be the vector space of all real polynomials of degree $\leq r$. The set $\mathcal{Z}$ being countable, we denote its elements by $\vec{z}_i, i \in \N$.
Let $V_k = \{\phi \in V: \phi(\vec{z}_i) = 0 \ \textrm{for} \ \textrm{all} \ i \leq k\}$.
Since $\dim V < \infty$, $V_k \supseteq V_{k+1}$, and $\mathcal{Z}$ is Zariski dense in $\R^m$, there exists $l$ such that $V_l = \{0\}$. In view of (\ref{2}), there exists $i_0 > 0$ such that
\begin{equation}
\mathrm{p}(\tau_i g\vec{z}_j) = \mathrm{p}(g\vec{z}_j)
\end{equation}
whenever $1 \leq j \leq l$ and $i \geq i_0$. Therefore
$$
\mathrm{p}(\tau_{i_0} g\vec{z}_j) = \mathrm{p}(g\vec{z}_j), 1 \leq j \leq l,
$$
proving that $\mathrm{p}(\tau_{i_0} g\vec{x}) - \mathrm{p}(g\vec{x}) \in V_l$.
So, $\mathrm{p}(\tau_{i_0}\vec{x}) = \mathrm{p}(\vec{x})$ and
$\tau_{i_0} \in L$. By (\ref{1}),
$$
\pi(g) = \tau_{i_0}^{-1}h_{i_0}\pi(e) \in L\pi(e).
$$
\qed

\subsection{Proof of Theorem \ref{thm0}}
If $f(\mathcal{Z})$ is discrete, it follows from Proposition \ref{connection1} that $H\pi(e)$ is closed. Hence Theorem \ref{thm0}  is an immediate consequence of the next proposition.

\begin{prop}
\label{auxiliary1} If $H\pi(e)$  is closed, then $f(\vec{x})$ is $\Q$-equivalent to a form in $n$ variables.
\end{prop}
{\bf Proof.} Suppose that $m >n$. Let $U_0$ be the unipotent radical of $H_0$. Note that
$$U_0 = \left\{\left(\begin{array}{cc} I_n&0\\ a&I_{m-n}\\ \end{array} \right): a \in \mathrm{M}_{m-n,n}(\R) \right\}.$$
Denote $W_0 = \{\vec{v} \in \R^m : h\vec{v} = \vec{v} \ \textrm{for} \ \textrm{all} \ h \in U_0\}$. Since $f_0 \in \Q[\vec{x}]$, $H_0$ is defined over $\Q$. Therefore $U_0$
is defined over $\Q$, and $U_0(\Q)$ is Zariski dense in $U_0$. Obviously, $W_0$ is also defined over $\Q$. In view of (\ref{stab}), $H_0$ is a semidirect product of a torus of dimension $n - 1$ and the real algebraic subgroup $H_{0u}$ generated by the unipotent elements of $H_0$.
Since $H = \sigma H_0 \sigma^{-1}$, we get that $U = \sigma U_0 \sigma^{-1}$ is the unipotent radical of $H$, $W = \{\vec{v} \in \R^m : h\vec{v} = \vec{v} \ \textrm{for} \ \textrm{all} \ h \in U\} = \sigma W_0$, and $H$ is also a semidirect product of a torus and the subgroup $H_u = \sigma {H_0}_u \sigma^{-1}$ generated by the unipotent elements of $H$.

By Ratner's theorem \cite{Ratner}, 
$\overline{H_u \pi(e)} = L\pi(e)$, where $L$ is a closed connected subgroup of $G$ containing $H_u$, and $L \cap \Gamma$ is a lattice in $L$. Since the orbit $H\pi(e)$ is closed, $L\pi(e) \subset H\pi(e)$. Therefore, $L$ is a subgroup of $H$ containing $H_u$. Since $H/H_u$ is a torus, $L/H_u$ consists of semi-simple elements. By \cite[Theorem 1]{Toma1} any  quotient of $L$ by its proper normal algebraic subgroup contains a unipotent element different from the identity. Therefore $L/H_u$ is trivial, i.e., $L = H_u$.
Since $H_u \cap \Gamma$ is a lattice in $H$ and $H_u \cap \Gamma$ is Zariski dense in $H_u$, we get that $H_u$ is defined over $\Q$.
  So, $U$ is also defined over $\Q$, implying that $U(\Q)$ is Zariski dense in $U$.  Hence $W = \{\vec{v} \in \R^m : h\vec{v} = \vec{v} \ \textrm{for} \ \textrm{all} \ h \in U(\Q)\}$, and $W$ is defined over $\Q$. But $\dim W_0 = \dim W$, and both $W_0$ and $W$
are defined over $\Q$. Therefore
  there exists a $\sigma_1 \in G(\Q)$ such that $\sigma_1 W = \sigma_1\sigma W_0 = W_0$. Replacing $f$ by its $\Q$-equivalent form $\sigma_1 f$ we may (and will) assume that $W = W_0$. It remains to prove that $f$ depends on $n$ variables. If $\mathcal{N}_G(U_0)$ is the normaliser of $U_0$ in $G$, it is easy to see that $\mathcal{N}_G(U_0) = \{g \in G: g W_0 = W_0\}$.
Since $W = W_0$, we have $\mathcal{N}_G(U_0) = \mathcal{N}_G(U)$. But $U_0$  coincides with the unipotent radical of $\mathcal{N}_G(U_0)$, and $U$ coincides with the unipotent radical of $\mathcal{N}_G(U)$. Therefore $U = U_0$.
Let $\vec{v} = (v_1, \dots, v_m) \in \R^m$. Denote $\vec{v}_1 = (v_1, \dots, v_n, 0, \dots, 0) \in \R^m$ and  $\vec{v}_2 = \vec{v} - \vec{v}_1$. Put $\vec{v'}_1 = (v_1, \dots, v_n) \in \R^n$ and $\vec{v'}_2 = (v_{n + 1}, \dots, v_m) \in \R^{m - n}$.
Let  $M$ be the set of all $\vec{v}$ with non-zero $\vec{v'}_1$. Clearly, for every $\vec{v} \in M$
there exists $a \in \mathrm{M}_{m-n,n}(\R)$ such that $\vec{v'}_2 + a\vec{v'}_1 = \vec{0}$.
 If $u = \left(\begin{array}{cc} I_n&0\\ a&I_{m-n}\\ \end{array} \right)$, then $u\vec{v} = \vec{v}_1$. Since $u \in H$ and $uf = f$, we get that $f(\vec{v}) = f(\vec{v}_1)$ whenever $\vec{v} \in M$. But $M$ is Zariski dense in $\R^m$. Therefore
$$
f(x_1, \dots, x_m) = f(x_1, \dots, x_n, 0, \dots, 0),
$$
proving the proposition. \qed

\section{Proof of Theorem \ref{thm++}} \label{reduction}
From now on we suppose that $m = n$.
\subsection{Compactness criteria.}

The next two propositions provide compactness criteria.
\begin{prop}
\label{auxiliary2}  $H\pi(e)$ is compact if and only if $f(\mathcal{Z})$ is discrete and $f$ does not represent $0$ over $\Q$ non-trivially.
\end{prop}

{\bf Proof.} Let $H\pi(e)$ be compact. Then $H = \mathcal{H} \Delta$, where $\mathcal{H}$ is a compact subset of $H$ and $\Delta$ is defined by (\ref{Delta}).
Let $\vec{z}_i \in \mathcal{Z}
$ and $f(\vec{z}_i) \rightarrow a$. Since $f = \sigma f_0$, we have $f_0(\sigma^{-1}\vec{z}_i) \rightarrow a$. In view of (\ref{form}), there exist $h_i \in H_0$ such that the sequences $h_i \sigma^{-1} \vec{z}_i$ and $\sigma h_i \sigma^{-1} \vec{z}_i$ are bounded.
Since  $\sigma h_i \sigma^{-1} \in H$, we get $\sigma h_i \sigma^{-1} = \tau_i \delta_i$, where $\tau_i \in \mathcal{H}$ and $\delta_i \in \Delta$. Then the sequence  $\delta_i\vec{z}_i \in \mathcal{Z}$ is bounded, and, passing to a subsequence, we may (and will) assume that $\delta_i\vec{z}_i = \delta_1\vec{z}_1$ for all $i$. We get
$$f(\vec{z}_i) = f(\sigma h_i \sigma^{-1} \vec{z}_i) = f(\tau_i\delta_i\vec{z}_i) = f(\delta_i\vec{z}_i) = f(\vec{z}_1) = a,$$
proving the discreteness of $f(\mathcal{Z})$.

Let $\vec{z} \in {\Z}^n$ be such that $f(\vec{z}) = 0$. Put $\vec{z}_i = i\vec{z}$, $i \in \N$. In a similar way as above, we prove the existence of $h_i \in H_0$, $\tau_i \in \mathcal{H}$, and $\delta_i \in \Delta$ such that $h_i\sigma^{-1}\vec{z}_i \rightarrow \vec{0}$ and
$$
\sigma h_i\sigma^{-1}\vec{z}_i = \tau_i \delta_i \vec{z}_i \rightarrow \vec{0}.
$$
By the compactness of $\mathcal{H}$ and the discreteness of ${\Z}^n$, we get $\vec{z} = \vec{0}$ proving that
$f$ does not represent $0$ over $\Q$ non-trivially.

Now, assume that $f(\mathcal{Z})$ is discrete and $f$ does not represent $0$ over $\Q$ non-trivially. It follows from Proposition \ref{connection1} that $H\pi(e)$ is closed. Since $\mathcal{Z}$ is a union of level
 sets $\mathcal{Z}_d, d \in \N$, it is enough to consider the case when $\mathcal{Z} = \mathcal{Z}_d$. If, on the contrary,  $H\pi(e)$ is not compact, Mahler's compactness criterion
implies the existence of a sequence $g_i \in H$, and of a sequence of primitive vectors $\vec{z}_{i} \in \Z^n$  such that $g_i \vec{z}_{i} \rightarrow \vec{0}$. Since $d\vec{z}_{i} \in \mathcal{Z}_d$ and $f(\mathcal{Z}_d)$ is discrete, we get for large $i$ that
$$f(d\vec{z}_{i}) = f(\vec{z}_{i}) = 0,$$
contrary to the assumption that $f$ does not represent $0$ over $\Q$ non-trivially.
\qed

\begin{prop}
\label{anisotropic1} Let $\Delta_0$ be a subgroup of $\Delta$ of finite index without torsion,
and let
$\mathcal{A}$ be the subspace of the vector space $\mathrm{M}_n(\Q)$ generated by
$\Delta_0$. The orbit $H\pi(e)$ is compact if and only if one of the following holds:
\begin{enumerate}
  \item $\mathcal{A}$ is a number field of degree $n$, which is not a $CM$ field, and $H$ is a $\Q$-anisotropic $\Q$-torus such that $H(\Q) = \mathcal{A} \cap G$;
  \item $n$ is even, $\mathcal{A}$ is a totally real number field of degree $\frac{n}{2}$, and $H_d$ is a $\Q$-anisotropic $\Q$-torus with $H_d(\Q) = \mathcal{A} \cap G$.
\end{enumerate}
\end{prop}

{\bf Proof.} If (1) or (2) holds, then $H\pi(e)$ is compact by the definition of $H$ and the properties of the $\Q$-anisotropic tori (see subsection $2.3$).

Let $H\pi(e)$ be compact.
As a subgroup of $H$, $\Delta_0$ is diagonalizable over $\C$, which implies that $\mathcal{A}$ is a commutative semisimple $\Q$-algebra. By a theorem of Wedderburn \cite[Theorem 1.4.4]{Herstein} $\mathcal{A} = F_1 \oplus \dots \oplus F_r$, where $F_i$ are simple commutative artinian rings, and by a theorem of Wedderburn-Artin \cite[Theorem 2.1.6]{Herstein} every $F_i$ is a number field. Suppose that $r \geq 2$, and let $d_1 \in F_1^*$ and $d_2 \in F_2^*$.
Then $\R^n = \ker (d_1) \oplus \mathrm{im} (d_1)$, $\mathrm{im}(d_1) \neq \{0\}$, and $\mathrm{im}(d_2) \neq \{0\}$. Since $d_1\cdot d_2 = 0$, $\mathrm{im}(d_2) \subset \ker (d_1)$, implying that $\ker (d_1) \neq \{0\}$.
Moreover, $d_1$ commutes element-wise with $H$, and $d_1 \in \mathrm{M}_n(\Q)$. Hence
$\ker (d_1)$ and $\mathrm{im} (d_1)$ are proper $H$-invariant $\Q$-subspaces of $\R^n$. Let $\phi: H \rightarrow \R^{*2}, h \mapsto (\det(h|_{\ker (d_1)}), \det(h|_{\mathrm{im}(d_1)}))$, and let $H'$ be the Zariski closure of $\Delta_0$ in $H$.
Since $\Delta_0$ consists of matrices with integer coefficients and the restriction of $\phi$ to $H'$ is a $\Q$-rational homomorphism, $\phi(\Delta_0) \subseteq \{\pm 1\} \times \{\pm 1\}$.
Moreover, $H/\Delta_0$ is compact and $H$ is Zariski connected. Therefore, $\mathrm{im} (\phi)$ is a compact Zariski connected subgroup of $\R^* \times \R^*$, implying that $\mathrm{im} (\phi) = \{(1,1)\}$. So, $\dim H < n-1$, contradicting the maximality of the torus $H$. We have proved that $\mathcal{A}$ is a subfield of $\mathrm{M}_n(\Q)$ containing the scalar matrices in $\mathrm{M}_n(\Q)$.
Let $q = \dim_{\Q} \mathcal{A}$. It follows from Proposition \ref{fields} that $q$ divides $n$, i.e., $n = q \cdot l$.

It is easy to see that $H_0$ (as in subsection 2.2) coincides with its centraliser $\mathcal{C}_G(H_0)$ in $G$.
Thus, $H = \mathcal{C}_G(H)$. Every element of $\mathcal{A} \cap G$ is a linear combination of elements from $\Delta_0$. Since $\Delta_0 \subset H$, we get that  $\mathcal{A} \cap G \subset H(\Q)$, implying $\Delta_0 \subset \mathcal{A} \cap \Gamma \subset \Delta$. Note that $\mathcal{A} \cap \mathrm{M}_n(\Z)$ is a sub-ring of finite index of the ring of integers $\mathcal{O}_{\mathcal{A}}$ of the number field $\mathcal{A}$. It follows that $\mathcal{A} \cap \Gamma$ and, therefore, $\Delta_0$ has finite index in the unit group $\mathcal{O}_{\mathcal{A}}^*$. Let $[\mathcal{A} : \Q] = 2s_1 + t_1$, where $s_1$ is the number of pairs of complex conjugated embeddings of $\mathcal{A}$ and $t_1$ is the number of real embeddings of $\mathcal{A}$. By Dirichlet's unit theorem,
the rank of $\Delta_0$ is equal to $s_1 + t_1 - 1$. On the other hand, $\Delta_0$ is of rank $s+t-1$ because $\Delta$ is a lattice of $H$ and $H \cong \mathrm{SO}_2(\R)^s \times (\R_{>0})^{s+t-1}$.
Consequently,
\begin{equation}\label{imp1}
  s + t = s_1 + t_1,
\end{equation}
and
\begin{equation}\label{imp2}
n = 2s + t = l[\mathcal{A} : \Q] = l(2s_1 + t_1) = l(s + t + s_1).
\end{equation}
Note that $l \leq 2$, because if $l \geq 3$, then $2s + t < l(s + t + s_1)$, contradicting (\ref{imp2}). If $l = 1$, then, in view of  (\ref{imp1}) and (\ref{imp2}), $s = s_1$, $t = t_1$, and $[\mathcal{A} : \Q] = n$. In this case, $\mathcal{A}$ is a maximal subfield of $\mathrm{M}_n(\Q)$ and, therefore, $H(\Q) = \mathcal{A} \cap G$ and $H$ is a $\Q$-anisotropic torus. By the choice of $\Delta_0$, the field $\mathcal{A}$ is generated by any subgroup of finite index of $\Delta_0$. Proposition \ref{CM-field} implies that $\mathcal{A}$ is not a $CM$ field. We have proved part $(1)$ in the formulation of the proposition.
Similarly, if $l = 2$
it follows from (\ref{imp1}) and (\ref{imp2}) that $t + 2s_1 = 0$. Hence, $s_1 = t = 0$ and $n = 2s$, proving that the torus $H_d$ is $\Q$-anisotropic and $H_d(Q) = \mathcal{A} \cap \Gamma$, where $\mathcal{A}$ is a totally real subfield of $\mathrm{M}_n(\Q)$ of degree $\frac{n}{2}$. This completes the proof of $(2)$ in the formulation of the proposition.
\qed

\subsection{Characterization of $f$ in terms of $H$.}

We have the following.
\begin{prop}
\label{anisotropic2}
\begin{enumerate}
  \item $f$ is a norm form if and only if $H$ is a $\Q$-anisotropic $\Q$-torus. In this case, $f$ corresponds to a maximal subfield of $\mathrm{M}_n(\Q)$ generated by $H(\Q)$.
  \item $f$ is a quasi-norm form if and only if $n = 2(\dim H_d + 1)$ and $H_d$ is a $\Q$-anisotropic $\Q$-torus. In this case, $f$ corresponds to a totally real number subfield of $\mathrm{M}_n(\Q)$ of degree $\frac{n}{2}$ generated by any subgroup of $\Delta$ of finite index without torsion.
\end{enumerate}
\end{prop}

{\bf Proof.} In order to prove (1), first suppose that $f$ is a norm form $\Q$-equivalent to $\mathfrak{N}_K$.
Then $H$ is a $\Q$-torus, and by Dirichlet's unit theorem $H(\Z)$ is co-compact in $H$, meaning that $H$ is $\Q$-anisotropic. Conversely, suppose that $H$ is $\Q$-anisotropic.  Let $\mathcal{K}$ be the sub-space of $\mathrm{M}_n(\Q)$ generated by $H(\Q)$.
Since $H = \mathcal{C}_G(H)$, $\mathcal{K} \cap G \subset H(\Q)$, implying
$\mathcal{K} \cap G = H(\Q)$. By Wederburn's theorem, $\mathcal{K} = \mathcal{K}_1 \oplus \dots \oplus \mathcal{K}_r$, where $\mathcal{K}_i$ are fields. If $r \geq 2$, then the map $H(\Q) \rightarrow \Q^*$, which sends every $g \in H(\Q)$ to the algebraic norm of the projection of $g$ in $\mathcal{K}_1$, gives rise to a non-trivial $\Q$-character of the torus $H$. This contradicts the assumption that $H$ is $\Q$-anisotropic.
Therefore $\mathcal{K}$ is a maximal subfield of $\mathrm{M}_n(\Q)$ and $f$ is a norm form corresponding to $\mathcal{K}$, which completes the proof of (1).

Let us prove (2). Suppose that $f$ is a quasi-norm form corresponding to a totally real number field $F$ of degree $s = \frac{n}{2}$. Let $F^1$ be the subgroup of $F^*$ of elements with algebraic norm $1$, that is, $$F^1 = \{x \in F: \theta_1(x) \dots \theta_s(x) = 1\}.$$ Then $F^1$ is
the group of $\Q$-points of a $\Q$-anisotropic torus of dimension $s-1$, which splits over $\R$ and fixes $f$, i.e., $F^1 \subset H_d$.
Since the rank of the group of units $\mathcal{O}_F^*$ is equal to the dimension of $H_d$, we have that $H_d \cap \Gamma$ is a lattice of $H_d$ and $H_d$ is $\Q$-anisotropic of dimension $s - 1$.

Suppose that $n = 2(\dim H_d + 1)$ and $H_d$ is $\Q$-anisotropic. Then $H\pi(e)$ is compact, and, in view of Proposition \ref{anisotropic1}, $H_d(\Q) = \mathcal{A} \cap G$, where $\mathcal{A}$ is a totally real subfield of $\mathrm{M}_n(\Q)$ of degree $s = \frac{n}{2}$ generated by a subgroup  of $\Delta$ of finite index without torsion. Fix an injective homomorphism of $\mathcal{A}$ into $\mathrm{M}_s(\Q)$, and, for every
 $a \in \mathcal{A}$, denote by $[a]$ the image of $a$ in $\mathrm{M}_s(\Q)$ under this injection.
We get the following homomorphism of the field $\mathcal{A}$ in $\mathrm{M}_n(\Q)$: for every $a \in \mathcal{A}$ the image of
$a$ in $\mathrm{M}_n(\Q)$ is given by the matrix
\begin{equation}
\left(\begin{array}{cc} [a]&0\\0&[a]\\ \end{array}\right) \in \mathrm{M}_n(\Q).
\end{equation}
By the Noether-Scolem theorem \cite[Theorem 4.3.1]{Herstein}, the image of $\mathcal{A}$ in $\mathrm{M}_n(\Q)$ is conjugate over $\Q$ to $\mathcal{A}$ itself. Therefore, $\R^n = W_1 \oplus W_2$, where $W_i$ are $\mathcal{A}$-invariant $\Q$-vector subspaces of $\R^n$ such that $\mathcal{A}$ acts in the same way on each of them; in other words, there exists an isomorphism of $\mathcal{A}$-modules $\omega: W_1 \rightarrow W_2$. Fix a $\Q$-basis $v_1, \dots, v_s$ of $W_1$ and let $v_{s+1} = \omega(v_1), \dots, v_{n} = \omega(v_s)$.
Denote $\vec{x} = (\vec{x}_1, \vec{x}_2)$, where $\vec{x}_1 = (x_1, \dots, x_s)$ is the coordinate vector with respect to the basis $v_1, \dots, v_s$ and $\vec{x}_2 = (x_{s+1}, \dots, x_n)$ is the coordinate vector with respect to the basis $v_{s + 1}, \dots, v_n$. Since $\mathcal{A}$ is a totally real number field of degree $s$, the natural action of $H_d(\Q)$ on every $W_i$ is diagonalizable over $\mathcal{A}$. Therefore, if $W_i^\star$ is the space dual to $W_i$, there exist linear forms $l_1, \dots, l_s$ in $s$ variables and uniquely determined pairwise different real characters $\chi_1, \dots, \chi_s$ of the torus $H_d$ such that $$W_i^\star = \R l_1(\vec{x}_i) + \dots +\R l_s(\vec{x}_i)$$ and $g l_i = \chi_i(g)l_i$ for the natural action of $H_d$ on $W_i^\star$. Hence, $V_i^\star = \R l_i(\vec{x}_1) + \R l_i(\vec{x}_2)$, $i = 1, \dots, s$, are the weight subspaces for the action of $H_d$ on the dual space ${\R^n}^\star$ of $\R^n$.
Since $H_d$ is a subgroup of $\mathrm{SL}_n(\R)$ and $H_d$ is Zariski connected, we have that
\begin{equation}\label{char}
\underset{i = 1}{\overset{s}{\prod}}\chi_i(g) = 1,
\end{equation}
for all $ g \in H_d$.
On the other hand, $H_c \cong  \underset{s}{\underbrace{\mathrm{SO}_2(\R) \times \dots \times \mathrm{SO}_2(\R)}}$,
 which implies that ${\R^n}^\star = V_1' \oplus \dots \oplus V_s'$, where $V_i'$, are uniquely defined $2$-dimensional $H_c$-invariant subspaces.
But $H_c$ and $H_d$ commute elementwise, and $V_i'$ are $H_c$-irreducible. Therefore, $V_i'$ are $H_d$-invariant and $H_d$ acts by scalar matrices on every $V_i'$; that is, $V_i'$ are weight subspaces for the action of $H_d$ on ${\R^n}^\star$. After permutation of indices, we have $V_i^\star = V_i'$. It follows that for every $i$ there exists a positive definite quadratic form $q_i$ on the $2$-dimensional subspace $V_i^\star$ such that $q_i(l_i(\vec{x}_1), l_i(\vec{x}_2))$ is $H_c$-invariant. It follows from (\ref{char}) that the quasi-norm form defined by (\ref{def2}) is $H$-invariant. 
Proposition \ref{prop.inv} implies that $f$ is a quasi-norm form corresponding to the totally real field $\mathcal{A}$. \qed

\begin{remark}\label{form=quasi/form}
Suppose that $f$ is both a norm and a quasi-norm form. In view of Proposition \ref{anisotropic2}, as a norm form, $f$ corresponds to the maximal subfield $\mathcal{K}$ of $\mathrm{M}_n(\Q)$ generated by $H(\Q)$ and, as a quasi-norm form, $f$ corresponds to the totally real subfield $F$ of $\mathrm{M}_n(\Q)$ of degree $\frac{n}{2}$ generated by a subgroup of finite index of $\Delta$. Since $F \subsetneqq \mathcal{K}$, it follows from Proposition \ref{CM-field} that $\mathcal{K}$ is a $CM$ field.
\end{remark}

\subsection{Proof of Theorem \ref{thm++}} According to Proposition \ref{auxiliary2}, condition $(2)$ in the formulation of
 Theorem \ref{thm++} is equivalent to the compactness of $H\pi(e)$.  By Proposition \ref{anisotropic1},  $H\pi(e)$ is compact if and only if either $H$ is $\Q$-anisotropic or $H_d$ is $\Q$-anisotropic. Finally, the equivalence between $(1)$ and $(3)$ in the formulation of
 Theorem \ref{thm++}, as well as the assertions $(a)$ and $(b)$ of $(3)$, follow immediately from Proposition \ref{anisotropic2}. \qed

\section{Proof of Theorem \ref{cor2}}

\subsection{Definition of a class of non-totally real forms}
In view of Proposition \ref{number fields}, given $n = 2s + t \geq 3$ with $s > 0$, there exists a number field $K$ of degree $n$ which is not a $CM$ field if $t = 0$, that is, if $K$ is totally imaginary.
The norm form $\mathfrak{N}_K$ can be writhen as $\mathfrak{N}_K = \phi_s f_t$, where
$\phi_s$ (respectively, $f_t$) is a product of non-real (respectively, real) linear forms. We have $\phi_s = (\lambda_{1}^2 + \lambda_{2}^2) \cdots (\lambda_{2s-1}^2 + \lambda_{2s}^2)$, where $\lambda_{i}$ are linearly independent real linear forms on $\R^n$.
Let $q_1, \dots, q_s$ be definite quadratic forms in two variables  such that the totally imaginary form $f_s = q_1(\lambda_{1}, \lambda_{2}) \cdots q_s(\lambda_{2s-1}, \lambda_{2s})$ is {not} proportional to $\phi_s$, that is, at least one of the forms in two variables $q_i(y_1,y_2)$ is not proportional to $y_1^2 + y_2^2$. \textit{Put} $f = f_s \cdot f_t$. Further on, we denote by $H^f$ the Zariski identity component of the stabilizer of $f$ in $G$ and by $H^f_c$ and $H^f_d$ the maximal compact and the maximal $\R$-split sub-tori of $H^f$. (See subsection $2.2$.)

\subsection{Proof of Theorem \ref{cor2}(1)} We will prove that the form $f$ introduced in subsection $5.1$ is as needed.
Using Mahler's criterion, it is easy to see that ${H^f\pi(e)}$ is relatively compact if and only if $f$ does not represent $0$ over $\Q$ non-trivially and $f(\mathcal{Z})$ is discrete at $0$. In view of Theorem \ref{thm++}, $f$ is neither a norm nor a quasi-norm form if and only if ${H^f\pi(e)}$ is not compact. Therefore, in order to prove Theorem \ref{cor2}(1) it is enough to show that ${H^f\pi(e)}$ is relatively compact but not compact.

 Let $H'$ be the maximal torus fixing $\mathfrak{N}_K$. Then $H' = H'_c \cdot H'_d$, where $H'_c$ is the maximal compact sub-torus of $H'$ and $H'_d$ is the maximal $\R$-split sub-torus of $H'$. It follows from the definition of $f$ and the description of $H_d$ in subsection $2.2$ that $H'_d = H^f_d$. Since $H' \pi(e)$ is compact, $\overline{H^f\pi(e)} = H^f_c \overline{H'_d\pi(e)}$ is also compact, i.e., ${H^f\pi(e)}$ is relatively compact.
It follows from the compactness of $H' \pi(e)$ that there
exists a compact connected subgroup $C$ of $H'_c$ such that $\overline{H^f_d \pi(e)} = C H^f_d \pi(e)$. Suppose, to the contrary, that $H^f\pi(e)$ is compact. Then
$$
H^f\pi(e) = H^f_c \overline{H^f_d \pi(e)} = H^f_c C  H^f_d\pi(e).
$$
In view of the discreteness of $\Gamma$, there exists a neighborhood $W$ of identity in $G$ such that
\begin{equation}\label{cor1.4}
H^f_c C  H^f_d \cap W \subset H^f.
\end{equation}
Note that $H'_c = \mathrm{SO}(\psi'_1) \times \dots \times \mathrm{SO}(\psi'_s)$, where $\psi'_i$ is the restriction of $\lambda_{2i-1}^2 + \lambda_{2i}^2$ to the $2$-dimensional $H^f_d$-invariant subspace $V_i$ of $\R^n$ (see subsection $2.2$).
Since $K$ is not a $CM$ field, it follows from Proposition \ref{ort.proj.} that $p_i(C) = \mathrm{SO}(\psi'_i)$ for all $i$.
On the other hand, $H^f_c = \mathrm{SO}(\psi_1) \times \dots \times \mathrm{SO}(\psi_s)$, where $\psi_i$ is the restriction of $q_i(\lambda_{2i-1}, \lambda_{2i})$ to $V_i$. By definition, the forms $f_s$ and $\phi_s$
are not proportional. Therefore there exists $1 \leq i_0 \leq s$ such that $\psi'_{i_0} \notin \R \psi_{i_0}$.
Note that $$p_{i_0}(H^f_c C H^f_d) = \mathrm{SO}(\psi_{i_0}) \mathrm{SO}(\psi'_{i_0})Z_{i_0} \subset \mathrm{GL}(V_{i_0}),$$ where $Z_{i_0}$ is the center of $\mathrm{GL}(V_{i_0})$.
Since $\mathrm{SO}(\psi_{i_0}) \neq \mathrm{SO}(\psi'_{i_0})$, the subgroup of $\mathrm{GL}(V_{i_0})$ generated by $p_{i_0}(H^f_c C H^f_d \cap W)$ is not abelian, contradicting (\ref{cor1.4}). \qed

\subsection{Proof of Theorem \ref{cor2}(2)} Let $H$ be a maximal non-diagonalizable over $\R$ algebraic torus of $G$.
Then $n = 2s + t$, where $s = \dim H_c > 0$ and $t = \dim H_d$. By the first part of the theorem, there exists a form $f$ of degree $n = 2s + t$ where $s$ (respectively, $t$) is the number of pairs of complex-conjugate (respectively, the number of real) linear forms dividing $f$. With $H^f$ as above, $H^f\pi(e)$ is relatively compact but not compact. It is clear that $H$ and $H_f$ have the same $\R$-rank, equal to $s + t -1$. Therefore, $H^f = g^{-1}Hg$, $g \in G$, proving that $H\pi(g)$ is also relatively compact but not compact. \qed

\end{document}